\documentclass{article}
\usepackage{amsfonts}
\usepackage{amsmath}
\usepackage{amsmath}
\usepackage{amssymb}

\newtheorem{theorem}{Theorem}

\newtheorem{lemma}[theorem]{Lemma}

\begin{document}

\title{Branching processes in random environment with sibling
dependence\thanks{This work is supported by the Program of the Presidium of the Russian Academy of Sciences No 01 'Fundamental Mathematics and its Applications' under grant PRAS-18-01.} }
\author{Vatutin
V.A.\thanks{Steklov Mathematical Institute, 8, Gubkin str., 119991, Moscow, Russia;
e-mail: vatutin@mi.ras.ru}, Dyakonova
E.E.\thanks{Steklov Mathematical Institute, 8, Gubkin str., 119991, Moscow, Russia;
e-mail: elena@mi.ras.ru}}

\date{}
\maketitle

\begin{abstract}
We consider a population of particles with unit life length. Dying each
particle produces offspring whose size depends on the random environment
specifying the reproduction law of all particles of the given generation and
on the number of relatives of the particle. We study the asymptotic behavior
of the survival probability of the population up to a distant moment $n$
under some restrictions on the properties of the environment and family ties.
\end{abstract}

\section{Introduction and main results}

We consider a population of particles with unit life length. Dying each
particle produces offspring whose size depends on a random environment
specifying the reproduction law of all particles of the given generation and
on the number of relatives of the particle. It will be clear from the
description to follow that the model we consider includes a class of
Galton-Watson branching processes in random environment (BPRE's).

First we give an informal description of the model. We fix a positive
integer $N$ and, for each $i\in \left\{ 1,...,N\right\} $ specify on the set
of $i$-dimensional vectors $\left( k_{1},...,k_{i}\right) $ with positive
integer components $0\leq k_{j}\leq N$ a probability measure ~$P(i;\cdot )$
with \begin{equation*}
\sum_{k_{1}=0,...,k_{i}=0}^{N}P(i;\left( k_{1},k_{2},...,k_{i}\right) )=1,
\end{equation*}such that\begin{equation}
P(i;\left( k_{1},k_{2},...,k_{i}\right) )=P\left( i;\left( k_{\sigma \left(
1\right) },k_{\sigma \left( 2\right) },...,k_{\sigma \left( i\right)
}\right) \right)  \label{Permut}
\end{equation}%
for any transposition $\left( \sigma \left( 1\right) ,\sigma \left( 2\right)
,...,\sigma \left( i\right) \right) $ of elements of the set $\left\{
1,...,i\right\} $.

It follows from this assumption that all the marginal distributions of the
probability measure $P(i;\cdot )$ coincide.

We say that a tuple \begin{equation*}
\mathbf{P}_{N}=\left\{ P(i;\cdot ),i\in \left\{ 1,...,N\right\} \right\}
\end{equation*}%
of probability measures on $\left\{ 1,...,N\right\} $ constitute an
environment of order $N$. \ The set of all environments
$\mathcal{P}_{N}=\left\{ \mathbf{P}_{N}\right\} $\ of order $N$ equipped with the metric
generated by the distance of total variation of the respective components of
$\mathbf{P}_{N}$ is a Polish space. Therefore we can consider probability
measures on this space. Let $\mathcal{P}$ \ be such a probability measure. A
sequence \begin{equation}
\mathbf{P}_{N}^{(n)}=\left( P^{(n)}(i;\cdot ),i\in \left\{ 1,...,N\right\}
\right) ,n=0,1,2,...,  \label{envirom}
\end{equation}%
of elements of $\mathcal{P}_{N}$, which are selected at random and
independently according to the measure $\mathcal{P},$ is said to form a
random environment.

A detailed description of the restriction we impose on the properties of the
random environment will be given later.

Suppose that a random environment $\left\{ \mathbf{P}_{N}^{(n)},n\geq
0\right\} \ $\ is fixed. With the environment in hands we can describe the
development of the BPRE $\left\{ \zeta (n),n\geq 0\right\} $ with sibling
dependence as follows. The process is initiated at time $n=0$ by $\zeta
(0)\geq 1$ particles of generation zero. We call these particles a sibling
group. Particles of the initial generation have a unit life length and dying
produce children. All the direct descendants of a particle of the initial or
subsequent generations \ will be called siblings (thus, siblings have one
and the same parent-particle). Some particles may have no direct
descendants. In this case the respective set of siblings is empty. The total
number of particles in a sibling group is called \textit{the} \textit{type
of the sibling group }or simply\textit{\ the} \textit{type of siblings}. If
a sibling group is of type $i$, we agree to consider that each particle in
this group has type $i$ as well.

If $\zeta (0)=i,$ then the initial particles die at moment $n=1$ and produce
in total $K(i)=k_{1}+\cdot \cdot \cdot +k_{i}$ particles of the first
generation, where $k_{j}$ is the number of descendants of the $j$-th
particle from the initial sibling group. The distribution of the vector
$\left( k_{1},k_{2},...,k_{i}\right) $ is specified by the measure
$P^{(0)}(i;\cdot ).$ Thus, $i$ new sibling groups are produced, each of which
contains all descendants (and only them) of some particle of zero
generation. If $k_{j}\geq 1$ for some $j$, then all the particles from this
sibling group die at moment $n=2$ and, independently of the behavior of
particles from the other sibling groups and the prehistory of the process
give birth to $K(k_{j})=i_{1}+\cdots +i_{k_{j}}$ particles of the second
generation, where $i_{r}$ is the number of descendants of the $r$th particle
from the $j-$th sibling group of the first generation. Thus, we have $k_{j}$
new sibling groups consisting of the second generation particles. The
distribution of the vector $\left( i_{1},i_{2},...,i_{k_{j}}\right) $ is
specified by the measure $P^{(1)}(k_{j};\cdot ),$ and so on....

Thus, given the environment the sibling groups existing at moment $n\geq 1$
evolve independently of each other. However the interaction of siblings at
this moment is described by (random) measures $P^{\left( n\right) }(i;\cdot
),i=1,...,N,$ specifying the joint distribution of the number of direct
descendants of a type $i$ sibling group.

Let $\zeta (n)$ denote the number of particles in generation $n$ in such
BPRE with sibling dependence. The aim of this note is to investigate the
asymptotic behavior of the survival probability of the process as
$n\rightarrow \infty $ under different conditions on the properties of random
environment.

We would like to note that there are several papers studying the behavior of
the Galton-Watson branching processes with sibling dependencies evolving in
a constant environment. We mention here only papers \cite{Bro88} and
\cite{Olof96}, which are the most significant for us. However, as far as we know,
BPRE's with sibling dependencies have not been yet analyzed.

According to the condition (\ref{Permut}) the marginal distributions of the
measure $P^{\left( n\right) }(i;\cdot )$ coincide for any $i=1,...,N$ .
Therefore, for any $j=0,1,...,N$ we can correctly define the (random )
variable \begin{equation}
p_{ij}^{\left( n\right) }:=\sum_{k_{2}=0,...,k_{i}=0}^{N}P^{\left( n\right)
}(i;\left( j,k_{2},...,k_{i})\right) ,  \label{defprob1}
\end{equation}%
which is equal the probability of the event that a particle, belonging at
time $n$ to some sibling group of type $i,$ begets just $j$ children, i.e.
generates a type $j$ sibling group.

We associate with the random environment (\ref{envirom}) two sequences of
(random) vector-valued multivariate generating functions \begin{equation}
\mathbf{\Phi }^{(n)}(\mathbf{s})=\left( \Phi _{1}^{(n)}(\mathbf{s}),...,\Phi
_{N}^{(n)}(\mathbf{s})\right) ,n\geq 0,  \label{3aa}
\end{equation}and \begin{equation*}
\mathbf{F}^{(n)}(\mathbf{s})=\left(
F_{1}^{(n)}(\mathbf{s}),...,F_{N}^{(n)}(\mathbf{s})\right) ,n\geq 0,
\end{equation*}%
where $\mathbf{s}=\left( s_{1},...,s_{N}\right) \in \left[ 0,1\right] ^{N}$
and, for $i=1,\ldots ,N,$ \begin{equation}
\Phi _{i}^{(n)}(\mathbf{s})=P^{(n)}(i;(0,...,0))+\sum_{\substack{
k_{1}=0,...,k_{i}=0,  \\ k_{1}+...+k_{i}>0}}^{N}P^{\left( n\right)
}(i;\left( k_{1},k_{2},...,k_{i}\right) )s_{k_{1}}s_{k_{2}}\cdot \cdot \cdot
s_{k_{i}},  \label{DefPhi}
\end{equation}(we assume that $s_{0}\equiv 1$) and

\begin{equation}
F_{i}^{(n)}(\mathbf{s})=p_{i0}^{\left( n\right)
}+\sum_{j=1}^{N}p_{ij}^{\left( n\right) }s_{j}^{j},\,i=1,\ldots ,N.
\label{ddefF}
\end{equation}%
Thus, the component $\Phi _{i}^{(n)}(\mathbf{s})$ of the vector-valued
multivariate generating function $\mathbf{\Phi }^{(n)}(\mathbf{s})$
describes in detail the joint law of generating \textbf{\ sibling groups }at
time $n$ by \textbf{all} \textbf{representatives} of a type $i$ sibling
group while the component $F_{i}^{(n)}(\mathbf{s})$ of the vector-valued
multivariate generating function $\mathbf{F}^{(n)}(\mathbf{s})$ describes
the distribution law of\textbf{\ the number of children} at time $n$ by
\textbf{a representative} of a type $i$ sibling group.

Recall, that the size of any sibling group (i.e., the number of children of
any particle) in our settings does not exceed $N$. Of course, this
assumption is an essential restriction. However, this assumption is natural
in the framework of applications in theoretical biology (see, for example,
monograph \cite{VJHBOOK}). Note, that it is more difficult (if possible at
all) to evaluate in practice parameters associated with the generating
function $\Phi _{i}^{(n)}(\mathbf{s})$ than to evaluate parameters
associated with the generating function $F_{i}^{(n)}(\mathbf{s})$. For this
reason we formulate the statements of theorems describing the asymptotic
behavior of the survival probability of BPRE's with sibling dependence in
the terms of the vector-valued (random) generating function \begin{equation*}
\mathbf{F}(\mathbf{s})=\left( F_{1}(\mathbf{s}),...,F_{N}(\mathbf{s})\right)
\end{equation*}with components \begin{equation}
F_{i}(\mathbf{s})=p_{i0}+\sum_{j=1}^{N}p_{ij}s_{j}^{j},\,i=1,\ldots ,N,
\label{dob1}
\end{equation}%
having the same distribution as the functions specified by (\ref{ddefF}).

We need some notation for $N$-dimensional vectors and $N\times N$ matrices.
Let $\mathbf{e}_{i}$, $i=1,2,\ldots ,N,$ be the $N$-dimensional vector whose
$i$-th component is equal to 1 and others are zeroes. For vectors
$\mathbf{x}=\left( x_{1},\ldots ,x_{N}\right) $ and $\mathbf{y}=\left( y_{1},\ldots
,y_{N}\right) $ set \begin{equation*}
(\mathbf{x},\mathbf{y}):=\sum_{i=1}^{N}x_{i}y_{i},\quad
|\mathbf{x}|:=\sum_{i=1}^{N}|x_{i}|\quad \text{and}\quad
\mathbf{x}^{\mathbf{y}}:=\prod_{i=1}^{N}(x_{i})^{y_{i}}.
\end{equation*}%
For a $N\times N$ matrix $\mathbf{m}=\left( m\left( i,j\right) \right)
_{i,j=1}^{N}$ introduce its norm by the equality \begin{equation*}
|\mathbf{m}|=\sum_{i=1}^{N}\sum_{j=1}^{N}|m\left( i,j\right) |.
\end{equation*}

Set $\delta _{kl}$ for the Kroneker symbol and let $\mathbf{1}=\left(
1,...,1\right) \in \left[ 0,1\right] ^{N}.$

Basic restrictions we impose on the properties of the BPRE with sibling
dependence are related with the mean matrix \begin{equation*}
\mathbf{M}=\mathbf{M}(\mathbf{F})=\left( M\left( i,j\right) \right)
_{i,j=1}^{N}:=\left( jp_{ij}\right) _{i,j=1}^{N}=\left(
\frac{\mathbb{\partial }F_{i}\left( \mathbf{1}\right) }{\partial s_{j}}\right) _{i,j=1}^{N}
\end{equation*}and the Hessian matrices \begin{equation*}
\mathbf{B}_{i}=\mathbf{B}_{i}(\mathbf{F})=\left( B_{i}(k,l)\right)
_{k,l=1}^{N}:=\left( k(k-1)p_{ik}\delta _{kl}\right) _{k,l=1}^{N}=\left(
\frac{\mathbb{\partial }^{2}F_{i}\left( \mathbf{1}\right) }{\partial
s_{k}\partial s_{l}}\right) _{k,l=1}^{N}
\end{equation*}%
constructed by the vector-valued generating function
$\mathbf{F}(\mathbf{s})$. The set of matrices generate two important random variables
\begin{equation}
\mathcal{B}:=\sum_{i=1}^{p}\left\vert \mathbf{B}_{i}\right\vert ,\qquad
\mathcal{T}:=\frac{\mathcal{B}}{\left\vert \mathbf{M}\right\vert ^{2}}.
\label{DefEta2}
\end{equation}

We define the cone \begin{equation*}
\mathcal{C=}\left\{ \mathbf{x}=(x_{1},...,x_{N})\in \mathbb{R}^{N}:x_{i}\geq
0\text{ for each }i\in \left\{ 1,...,N\right\} \right\} ,
\end{equation*}the sphere \begin{equation*}
\mathbb{S}^{N-1}=\left\{ \mathbf{x}:\mathbf{x}\in \mathbb{R}^{N},\left\vert
\mathbf{x}\right\vert =1\right\}
\end{equation*}and the space $\mathbb{X}=\mathcal{C}\cap \mathbb{S}^{N-1}$.

To go further we need to attract the linear semigroup $S^{++}$ of $N\times N$
matrices all whose elements are non-negative.

Assume that the distribution of the random matrix $\mathbf{M}$ meets the
following restrictions:

\textbf{Condition H1}. There exists $\theta >0$ such that \begin{equation*}
\mathcal{E}\left[ \left\vert \mathbf{M}\right\vert ^{\theta }\right]
=\int_{S^{+}}\left\vert \mathbf{M}\right\vert ^{\theta }\mathcal{P}\left(
d\mathbf{M}\right) <\infty .
\end{equation*}

\textbf{Condition H2}. (Strong irreducibility). The support of $\mathcal{P}$
in $S^{++}$ acts strongly irreducibly on $\mathbb{R}^{N},$ i.e. no proper
finite union of subspaces of $\mathbb{R}^{N}$ is invariant with respect to
all elements of the multiplicative semi-group generated by the support of
$\mathcal{P}$.

\textbf{Condition H3}. Elements of the random matrix $\mathbf{M}=\left(
M\left( i,j\right) \right) _{i,j=1}^{N}$ are positive and there exists a
real positive number $\gamma >1$ such that $\mathcal{P}$-a.s.
\begin{equation*}
\frac{1}{\gamma }\leq \frac{M\left( i,j\right) }{M\left( k,l\right) }\leq
\gamma \end{equation*}for any $i,j,k,l\in \left\{ 1,...,N\right\} $.

For $n=0,1,2,\ldots $ introduce random matrices \begin{equation}
\mathbf{M}^{(n)}=\mathbf{M}^{(n)}(\mathbf{F}^{(n)})=\left( M^{\left(
n\right) }\left( i,j\right) \right) :=\left( jp_{ij}^{\left( n\right)
}\right) =\left( \frac{\mathbb{\partial }F_{i}^{(n)}\left( \mathbf{1}\right)
}{\partial s_{j}}\right) _{i,j=1}^{N},  \label{matrix_n}
\end{equation}\begin{equation*}
\mathbf{B}_{i}^{(n)}=\mathbf{B}_{i}^{(n)}(\mathbf{F}^{(n)})=\left(
B_{i}^{(n)}(k,l)\right) _{k,l=1}^{N}:=\left( k(k-1)p_{ik}^{(n)}\delta
_{kl}\right) _{k,l=1}^{N}=\left( \frac{\mathbb{\partial
}^{2}F_{i}^{(n)}\left( \mathbf{1}\right) }{\partial s_{k}\partial s_{l}}\right) _{k,l=1}^{N},
\end{equation*}and denote by \begin{equation}
\mathbf{R}^{(n)}=\left( R^{(n)}(i,j)\right)
_{i,j=1}^{N}:=\mathbf{M}^{(0)}\mathbf{M}^{(1)}\cdots \mathbf{M}^{(n)}  \label{ProdR}
\end{equation}%
the right product of the random mean matrices $\mathbf{M}^{(n)},n\geq 0$.

It is known (see \cite{FK1960}) that given \begin{equation}
\mathcal{E}\left[ \max (0,\log \left\vert \mathbf{M}\right\vert )\right]
<\infty ,  \label{Fursten}
\end{equation}the sequence \begin{equation*}
\frac{1}{n}\log \left\vert \mathbf{R}^{\left( n\right) }\right\vert
,n=1,2,...
\end{equation*}%
converges $\mathcal{P}$-a.s. as $n\rightarrow \infty $ to a limit
\begin{equation*}
\Lambda :=\lim_{n\rightarrow \infty }\frac{1}{n}\mathcal{E}\left[ \log
\left\vert \mathbf{R}^{\left( n\right) }\right\vert \right] ,
\end{equation*}called the upper Lyapunov exponent.

We need to impose two more conditions on $\mathcal{P}$.

\textbf{Condition H4}. The upper Lyapunov exponent $\Lambda $ of the
distribution generated by $\mathcal{P}$ on $S^{++}$ is equal to $0$.

\textbf{Condition H5}. There exists $\delta >0$ such that \begin{equation*}
\mathcal{P}\left( \mathbf{M}\in S^{++}:\,\log \left\vert
\mathbf{xM}\right\vert \geq \delta \text{ for any }\mathbf{x}\in \mathbb{X}\,\right) >0.
\end{equation*}

We now are ready to formulate the first main result of the paper.

\begin{theorem}
\label{theor1} Assume Conditions $\mathbf{H1-H5}$. If \begin{equation}
\sup_{\mathbf{x}\in \mathbb{X}}\mathcal{E}\left[ \frac{1}{\left\vert
\mathbf{xM}\right\vert }\right] <\infty   \label{ExponFinite}
\end{equation}%
and, for $\mathcal{T}$ given in (\ref{DefEta2}) and an $\varepsilon >0$
\begin{equation}
\mathcal{E}\left[ \mathcal{T}^{1+\varepsilon }\right] <\infty ,
\label{SecondFinite}
\end{equation}%
then, for any $i=1,2,\ldots ,N$ there exists a number $\beta _{i}\in \left(
0,\infty \right) $ such that \begin{equation}
\lim_{n\rightarrow \infty }\sqrt{n}\mathcal{P}\left( \medskip \zeta
(n)>0\mathbf{|}\zeta (0)=i\right) =\beta _{i}.  \label{SurvivalProbab}
\end{equation}
\end{theorem}

In the sequel we call a BPRE with sibling dependence \textit{critical}, if
$\Lambda =0$.

We now describe conditions under which the asymptotics of the survival
probability of $\medskip $the process $\zeta (n)$\ has a form different from
that stated in Theorem \ref{theor1}.

\textbf{Condition }$\mathbf{H6}$. There exists an $\varepsilon >0$ such that
\begin{equation*}
\mathcal{E}\left[ |\log \mathcal{T}|^{1+\varepsilon }|\mathbf{M}|\right]
<\infty .
\end{equation*}

Let \begin{equation*}
\Theta :=\left\{ \theta >0:\mathcal{E}\left[ |\mathbf{M|}^{\theta }\right]
<\infty \right\} .
\end{equation*}%
It is known (see, for example, \cite{FK1960}) that for any $\theta \in
\Theta $ the limit \begin{equation*}
\lambda \left( \theta \right) :=\lim_{n\rightarrow \infty }\left(
\mathcal{E}\left[ |\mathbf{R}^{\left( n\right) }|^{\theta }\right] \right)
^{1/n}<\infty \quad \end{equation*}is well defined. Set \begin{equation*}
\Lambda (\theta ):=\log \lambda (\theta ),\quad \theta \in \Theta .
\end{equation*}

\begin{theorem}
\label{theor2} Assume that conditions $\mathbf{H1}-\mathbf{H3\ }$and
$\mathbf{H6}$ are valid, the point $\theta =1$ belongs to the interior of the
set $\Theta $ and $\Lambda ^{\prime }(1)<0$. Then, for any $i=1,\ldots ,N$

\begin{itemize}
\item[(a)] there exists a constant $C_{i}>0$ such that \begin{equation}
\mathcal{P}\left( \zeta (n)>0\Big|\zeta (0)=i\right) \sim C_{i}\lambda
^{n}(1),\quad n\rightarrow \infty ;  \label{P_asymp}
\end{equation}

\item[(b)] \begin{equation*}
\lim_{n\rightarrow \infty }\mathcal{E}\left[ s^{\zeta (n)}\,\Big|\zeta
(n)>0;\zeta (0)=i\right] =\Psi _{i}\left( s\right) ,\,s\in \lbrack 0,1),
\end{equation*}%
where $\Psi _{i}\left( s\right) $ is the probability generating function of
a proper distribution on $\mathbb{Z}_{+}$.
\end{itemize}
\end{theorem}

\subsection{Proofs of Theorems \protect\ref{theor1}--\protect\ref{theor2}}

The basic idea of proof of Theorems \ref{theor1} and \ref{theor2} is to
compare the BPRE $\left\{ \zeta (n),n\geq 0\right\} $ with sibling
dependence with another process, a macro process. This macro process
consists of sibling groups, to be called macro particles.

The type of a macro particle is the number of particles from the initial
BPRE with sibling dependence which belong to the sibling group constituting
the macro particle.

As we have mentioned, it will be convenient for us to assume that all the
particles, constituting the macro particle have the same type as the macro
particle. Thus, we assign each sibling group to one of $N$ possible types of
macro particles. This allows us to associate the BPRE $\left\{ \zeta
(n),n\geq 0\right\} $ with sibling dependence with the macro
process\begin{equation*}
\left\{ \mathbf{Z}\left( n\right) =\left(
Z_{1}(n),Z_{2}(n),...,Z_{N}(n)\right) ,\quad n\geq 0\right\} ,
\end{equation*}%
where $Z_{k}(n),k=1,2,...,N,$ is the number of macro particles of type $k$
in the $n$th generation of the macro process, i.e. $Z_{k}(n)$ is the number
of such sibling groups in the $n$th generation of $\zeta (n),$ each of which
is a sibling group of size $k,$ generated by a parent-particle belonging to
the $(n-1)$th generation.

Recall that an $i$ type macro particle existing at time $n$ in the macro
process generates offspring according to the probability generating function
$\Phi _{i}^{(n)}(\mathbf{s})$ specified in (\ref{DefPhi}) while the marginal
distributions for the number of direct descendants of a particle belonging
to a size $i$ sibling group are defined by (\ref{defprob1}).

Clearly, the macro process $\left\{ \mathbf{Z}\left( n\right) ,n\geq
0\right\} $ is an $N$-type Galton-Watson process in random environment.

The main difference between the processes $\left\{ \zeta (n),n\geq 0\right\}
$ and $\left\{ \mathbf{Z}\left( n\right) ,n\geq 0\right\} $ is easy to
explain: given the environment the individuals of the initial BPRE $\zeta
(n) $ with sibling dependence do not reproduce independently while macro
particles do, since the only dependencies are within the sibling groups.

It is not difficult to see that \begin{equation}
\zeta (n)=\sum_{k=1}^{N}kZ_{k}(n).  \label{base}
\end{equation}

We use the symbols $\mathcal{E}_{\mathbf{\Phi }},$
$\mathcal{D}_{\mathbf{\Phi }}$, and $\mathcal{P}_{\mathbf{\Phi }}$ for the expectations,
variances, and probabilities given the vector-valued probability generating
function $\mathbf{\Phi }$. Denote \begin{equation*}
\mathbf{M}_{macro}=\left( M_{macro}\left( i,j\right) \right)
_{i,j=1}^{N}:=\left( \mathcal{E}_{\mathbf{\Phi }}\left[ Z_{j}\left( 1\right)
;\mathbf{Z}(0)=\mathbf{e}_{i}\right] \right) _{i,j=1}^{N}
\end{equation*}the (random) mean matrix and by \begin{equation*}
\mathbf{B}_{i,macro}=\left( B_{i,macro}\left( j,k\right) \right)
_{j,k=1}^{N}:=\left( \mathcal{E}_{\mathbf{\Phi }}\left[ Z_{j}\left( 1\right)
\left( Z_{k}(1)-\delta _{jk}\right) ;\mathbf{Z}(0)=\mathbf{e}_{i}\right]
\right) _{j,k=1}^{N}
\end{equation*}the (random) Hessian matrices for the macro process. We put
\begin{equation*}
p_{2\left( jk\right) }:=P(2;\left( j,k,\right)
\end{equation*}and set, for $i\geq 3$ \begin{equation*}
p_{i\left( jk\right) }:=\sum_{(k_{3},...,k_{i})}P(i;\left(
j,k,k_{3},...,k_{i})\right) .
\end{equation*}

\begin{lemma}
\label{L1}The elements of the matrices $\mathbf{M}_{macro}$ and
$\mathbf{B}_{i,macro}$ are calculated by the formulas \begin{equation*}
\mathcal{E}_{\mathbf{\Phi }}\left[ Z_{j}\left( 1\right)
;\mathbf{Z}(0)=\mathbf{e}_{i}\right] =M_{macro}\left( i,j\right)
=ip_{ij}=\frac{i}{j}M\left( i,j\right) .
\end{equation*}\begin{equation*}
B_{i,macro}\left( j,k\right) =\mathcal{E}_{\mathbf{\Phi }}\left[ Z_{j}\left(
1\right) \left( Z_{k}(1)-\delta _{jk}\right)
;\mathbf{Z}(0)=\mathbf{e}_{i}\right] =i\left( i-1\right) p_{i\left( jk\right) }.
\end{equation*}
\end{lemma}

\textbf{Proof}. Let $\xi _{i,r}$ be the number of direct descendants of the
$r$th particle entering an $i$ type macro particle of generation zero and let
$I\left\{ \mathcal{C}\right\} $ be the indicator of the event \
$\mathcal{C}$.

Then\begin{eqnarray*}
\mathcal{E}_{\mathbf{\Phi }}\left[ Z_{j}\left( 1\right)
;\mathbf{Z}(0)=\mathbf{e}_{i}\right]  &=&\mathcal{E}_{\mathbf{\Phi }}\left[
\sum_{r=1}^{i}I\left\{ \xi _{i,r}=j\right\} \right]  \\
&=&i\mathcal{E}_{\mathbf{\Phi }}\left[ I\left\{ \xi _{i,1}=j\right\} \right]
=i\sum_{(k_{2},...,k_{i})}P(i;\left( j,k_{2},...,k_{i})\right) =ip_{ij}.
\end{eqnarray*}Recalling the condition (\ref{Permut}) we have for $r\neq
t$\begin{equation*}
\mathcal{E}_{\mathbf{\Phi }}\left[ I\left\{ \xi _{i,r}=j\right\} I\left\{
\xi _{i,t}=k\right\} \right] =\sum_{(k_{3},...,k_{i})}P(i;\left(
j,k,k_{3}...,k_{i})\right) =p_{i\left( jk\right) },
\end{equation*}and for $r=t$\begin{equation*}
\mathcal{E}\left[ I\left\{ \xi _{i,r}=j\right\} I\left\{ \xi
_{i,t}=k\right\} \right] =\delta _{jk}p_{ij}.
\end{equation*}

Therefore, the random Hessian matrices $\mathbf{B}_{i,macro}=\left(
B_{i,macro}\left( j,k\right) \right) _{j,k=1}^{N}$ have elements
\begin{eqnarray*}
B_{i,macro}\left( j,k\right)  &=&\frac{\partial ^{2}\Phi
_{i}(\mathbf{s})}{\partial s_{j}\partial s_{k}}=\mathcal{E}\left[ \left(
\sum_{r=1}^{i}I\left\{ \xi _{i,r}=j\right\} \right) \left(
\sum_{t=1}^{i}I\left\{ \xi _{i,t}=k\right\} -\delta _{jk}\right) \right]  \\
&=&\sum_{r\neq t}\mathcal{E}\left[ I\left\{ \xi _{i,r}=j\right\} I\left\{
\xi _{i,t}=k\right\} \right] +\sum_{r=1}^{i}\mathcal{E}\left[ I\left\{ \xi
_{i,r}=j;\xi _{i,r}=k\right\} \right]  \\
&&-\delta _{jk}\sum_{r=1}^{i}\mathcal{E}\left[ I\left\{ \xi _{i,r}=j\right\}
\right]  \\
&=&i(i-1)p_{i\left( jk\right) }+i\delta _{jk}p_{ij}-i\delta
_{jk}p_{ij}=i(i-1)p_{i\left( jk\right) }.
\end{eqnarray*}

The lemma is proved.

Let $\rho $ be the Perron root of the matrix $\mathbf{M}$ and let
$\mathbf{u}=\left( u_{1},...,u_{N}\right) $ be the right eigenvector of $\mathbf{M}$
corresponding to $\rho $ and satisfying the condition $\left\vert
\mathbf{u}\right\vert =u_{1}+...+u_{N}=1$.

\begin{lemma}
\label{L2}The value $\rho $ is the maximal in modulo eigenvalue of the
matrix $\mathbf{M}_{macro},$ and the right eigenvector $\mathbf{U}=\left(
U_{1},...,U_{N}\right) $ corresponding to $\rho $ has components
\begin{equation*}
U_{j}=\frac{ju_{j}}{\sum_{k=1}^{N}ku_{k}},j=1,...,N.
\end{equation*}
\end{lemma}

\textbf{Proof}. If \begin{equation*}
\sum_{j=1}^{N}M(i,j)u_{j}=\sum_{j=1}^{N}jp_{ij}u_{j}=\rho u_{i},
\end{equation*}then\begin{eqnarray*}
\sum_{j=1}^{N}M_{macro}(i,j)U_{j} &=&\sum_{j=1}^{N}\frac{i}{j}M\left(
i,j\right) \frac{ju_{j}}{\sum_{k=1}^{N}ku_{k}} \\
&=&\frac{i}{\sum_{k=1}^{N}ku_{k}}\sum_{j=1}^{N}M\left( i,j\right)
u_{j}=\frac{i}{\sum_{k=1}^{N}ku_{k}}\rho u_{i}=\rho U_{i}.
\end{eqnarray*}%
Thus, $\mathbf{U}$ is the right eigenvector corresponding to $\rho $.

Let us show that $\rho $ is the maximal in modulo eigenvalue of the matrix
$\mathbf{M}_{macro}$. Indeed, assume that there exists a $\rho ^{\ast }>\rho $
and the respective eigenvector \ $\mathbf{U}^{\ast }=\left( U_{1}^{\ast
},...,U_{N}^{\ast }\right) $ with strictly positive components such that
\begin{equation*}
\sum_{j=1}^{N}M_{macro}(i,j)U_{j}^{\ast }=\rho ^{\ast }U_{i}^{\ast }
\end{equation*}for all $i=1,2,...,N.$ Therefore,\begin{equation*}
\sum_{j=1}^{N}\frac{1}{j}M\left( i,j\right) U_{j}^{\ast }=\rho ^{\ast
}\frac{U_{i}^{\ast }}{i}.
\end{equation*}%
Hence, setting $U_{i}^{\ast \ast }:=i^{-1}U_{i}^{\ast },i=1,2,...,N,$ we see
that\begin{equation*}
\sum_{j=1}^{N}M\left( i,j\right) U_{j}^{\ast \ast }=\rho ^{\ast }U_{i}^{\ast
\ast }.
\end{equation*}%
Thus, $\rho ^{\ast }>\rho $ is an eigenvalue of $\mathbf{M}$. This
contradicts to the fact that $\rho $ is the maximal in modulo eigenvalue of
the matrix $\mathbf{M}$.

The lemma is proved.

\textbf{Proof of Theorem \ref{theor1}.} We consider, along with the macro
process, an auxiliary $N$-type branching process in random environment, the
so-called individual process. The reproduction of $i$-type particles at
moment $n$ in the new process is specified by the probability generating
function \begin{equation*}
F_{i}^{\left( n\right) }(s_{1},...,s_{N})=p_{i0}^{\left( n\right)
}+\sum_{j=1}^{N}p_{ij}^{\left( n\right) }s_{j}^{j},
\end{equation*}which is the same as in (\ref{ddefF}).

It follows from Lemma \ref{L1} that the mean matrix for the reproduction law
of the particles of the auxiliary process at moment $n$ has the form
\begin{equation*}
\mathbf{M}^{\left( n\right) }=\left( M^{\left( n\right) }\left( i,j\right)
\right) _{i,j=1}^{N}=\left( jp_{ij}^{\left( n\right) }\right) .
\end{equation*}

Thus, \begin{equation}
M_{macro}^{(n)}\left( i,j\right) =ip_{ij}^{(n)
}=\frac{i}{j}jp_{ij}^{(n)}=\frac{i}{j}M^{(n)}\left( i,j\right) .  \label{seee}
\end{equation}
Set \begin{equation*}
\mathbf{R}_{macro}^{(n)}=\left( R_{macro}^{(n)}(i,j)\right)
_{i,j=1}^{N}:=\mathbf{M}_{macro}^{(0)
}\mathbf{M}_{macro}^{(1)}\cdots\mathbf{M}_{macro}^{(n)}.
\end{equation*}

It is easy to see that the elements of the matrix $\mathbf{R}_{macro}^{(1)}$
satisfy the relation \begin{eqnarray*}
R_{macro}^{(1)}\left( i,j\right)
&=&\sum\limits_{k=1}^{N}\frac{i}{k}M^{(0)}\left( i,k\right) \frac{k}{j}M^{(1)}\left( k,j\right) \\
&=&\frac{i}{j}\sum\limits_{k=1}^{N}M^{(0)}\left( i,k\right) M^{(1)}\left(
k,j\right) =\frac{i}{j}R^{(1)}\left( i,j\right) .
\end{eqnarray*}

Hence we conclude by induction that the elements $R_{macro}^{\left( n\right)
}(i,j)$ of the matrix $\mathbf{R}_{macro}^{\left( n\right) }$ have the
form\begin{equation}
R_{macro}^{\left( n\right) }(i,j)=\frac{i}{j}R^{\left( n\right) }(i,j).
\label{macR}
\end{equation}Condition $\mathbf{H4}$ gives\begin{equation*}
\Lambda _{macro}:=\lim_{n\rightarrow \infty }\frac{1}{n}\mathcal{E}\left[
\log \left\vert \mathbf{R}_{macro}^{\left( n\right) }\right\vert \right]
=\Lambda =0.
\end{equation*}%
Using (\ref{macR}) it is easy to show that under the assumptions of Theorem
\ref{theor1} the macro process $\left\{ \mathbf{Z}\left( n\right) ,n\geq
0\right\} $ satisfies all the conditions of Theorem 1 in \cite{VD2016},
according to which \begin{equation*}
\mathcal{P}(\mathbf{Z}(n)\neq
\mathbf{0}\,|\,\mathbf{Z}(0)=\mathbf{e}_{i})\sim \frac{D_{i}}{\sqrt{n}}
\end{equation*}%
as $n\rightarrow \infty $. Since $\left\{ \mathbf{Z}(n)\neq
\mathbf{0}\right\} \Leftrightarrow \left\{ \sum_{k=1}^{N}kZ_{k}(n)>0\right\} ,$ we
conclude by (\ref{base}) that, as $n\rightarrow \infty $ \begin{equation*}
\mathcal{P}\left( \zeta (n)>0\,|\,\zeta (0)=i\right) =\mathcal{P}\left(
\sum_{k=1}^{N}kZ_{k}(n)>0|\mathbf{Z}(0)\mathbf{=e}_{i}\right) \sim
\frac{D_{i}}{\sqrt{n}}.
\end{equation*}Theorem \ref{theor1} is proved.

\bigskip \textbf{Proof of Theorem \ref{theor2}.} Based on the reasonings
used earlier in the proof of Theorem \ref{theor1} it is easy to check by
(\ref{macR}) that one may apply Theorem 1 in \cite{VW2018}, proved for the
so-called strongly subcritical multitype BPRE, to the branching macro
process $\left\{ \mathbf{Z}(n),n\geq 0\right\} $. This fact and the
equaliuty (\ref{base}) give the needed statement

\section{Limit distribution of the number of particles in the critical
BPRE's with sibling dependence}

To describe the limiting behavior of the distribution of the number of
particles in a critical BPRE with sibling dependence we need impose stronger
restrictions than in Theorem \ref{theor1}.

Given a random vector-valued generating function $\mathbf{\Phi
}(\mathbf{s})=\left( \Phi _{1}(\mathbf{s}),...,\Phi _{N}(\mathbf{s})\right)
=\mathbf{\Phi }^{\left( 0\right) }(\mathbf{s})$ from (\ref{3aa}) and (\ref{DefPhi})
we introduce a random vector $\left( \eta _{i1},...,\eta _{iN}\right)
,i=1,...,N,$ whose distribution is specified by the generating function
$\Phi _{i}(\mathbf{s})=\mathcal{E}s_{1}^{\eta _{i1}}...s_{N}^{\eta _{iN}}.$
Recall that $\mathcal{D}_{\mathbf{\Phi }}$ is the symbols for the variance
given the vector-valued probability generating function $\mathbf{\Phi }$ and
$\rho $ is the Perron root of the mean matrix $\mathbf{M}$, i.e. it is its
maximal in absolute value eigenvalue.

Introduce the random variables

\begin{equation*}
\Delta :=\max_{i,j}\mathcal{D}_{\mathbf{\Phi }}\left[ \eta _{ij}\right] .
\end{equation*}

We need the following conditions.

\textbf{Condition A1}. The mean matrices $\mathbf{M}$ and
$\mathbf{M}^{\left( n\right) },n\geq 0,$ are positive and with probability 1 have a
\textit{common non-random right eigenvector}
$\mathbf{u}=(u_{1},...,u_{p})^{\prime },\,\left\vert \mathbf{u}\right\vert =1$, with positive components
corresponding to their Perron roots $\rho $ and $\rho ^{(n)}$.

\textbf{Condition A2.} The following relation is valid \begin{equation*}
\mathcal{P}\left( \max_{1\leq i\leq N}\left( p_{i0}+p_{i1}\right)
<1\,\right) =1,
\end{equation*}where $p_{i0}$ and $p_{i1}$ are from (\ref{dob1}).

\textbf{Condition A3.} The distribution of the random variable $X:=\log \rho
$ belongs without centering to the domain of attraction of some stable law
$T $ with index $\alpha \in (0,2]$. The limit law $T$ is not one-sided law,
that is, $0<T\left( R^{+}\right) <1.$\

Note that according to Lemma \ref{L2} the Perron root of the mean matrix
$\mathbf{M}_{macro}$ is equal to the Perron root of $\mathbf{M}$.

\bigskip \textbf{Condition A4.} With probability 1 \begin{equation}
\Delta <\infty   \label{AsinfinB4}
\end{equation}and there exists $\varepsilon >0$ such that \begin{equation}
\mathcal{E}[\log ^{+}\rho ^{-2}\Delta ]^{\alpha +\varepsilon }<\infty ,\text{
\ }  \label{Ass22B4}
\end{equation}where $\alpha $ is from Condition A3.

Recall that the meander of a strictly stable Levy process is a strictly
stable Levy process conditioned to stay positive on the time interval $(0,1]$
(see \cite{Chau97} and \cite{Chau06}) for more detail.

\begin{theorem}
\label{theor3} Let Conditions A1-A4 be valid. Then there exists a slowly
varying at infinity sequence $l\left( 0\right) ,l\left( 1\right) ,l\left(
2\right) ,...$ such that, for any $i=1,\ldots ,N$ as $n\rightarrow \infty
$\begin{equation}
\lim_{n\rightarrow \infty }n^{-1/\alpha }l\left( n\right) \mathcal{P}\left(
\medskip \zeta (n)>0\,|\,\zeta (0)=i\right) =\beta _{i}\in \left( 0,\infty
\right) ,  \label{pria}
\end{equation}where $\alpha $ is from Condition A3, and \begin{equation}
\mathcal{L}\left( \left\{ n^{-1/\alpha }l\left( n\right) \log \zeta
([nt]),0\leq t\leq 1\right\} \,|\,\zeta (n)>0,\zeta (0)=i\right) \implies
\mathcal{L}\left( L^{+}\right) ,  \label{mash}
\end{equation}%
where $[x]$ denotes the integer part of the number $x$ and $L^{+}$ is the
meander of a strictly stable Levy process with index $\alpha .$
\end{theorem}

Here and what follows the symbol $\implies $ stands for the weak convergence
with respect to Skorokhod topology in the space $D[0,1]$ of cadlag functions
on the unit interval.

\textbf{Proof of Theorem \ref{theor3}. } We know that according to Lemma
\ref{L2} the Perron root of the mean matrix $\mathbf{M}_{macro}$ is equal to the
Perron root of $\mathbf{M}$. This and Condition A3 show that the macro
process constructed by the initial BPRE with sibling dependence satisfies
all the conditions of \ Lemma 9 and Theorem 1 \ in \cite{Dyak15}. Therefore,
there exists a slowly varying at infinity sequence $l\left( 0\right)
,l\left( 1\right) ,l\left( 2\right) ,...$ such that, for each $i=1,\ldots ,N$
as $n\rightarrow \infty $\begin{equation}
\lim_{n\rightarrow \infty }n^{-1/\alpha }l\left( n\right) \mathcal{P}\left(
\medskip \mathbf{Z}(n)>0\,|\,\mathbf{Z}(0)=\mathbf{e}_{i}\right) =\beta
_{i}\in \left( 0,\infty \right) ,  \label{iras}
\end{equation}where $\alpha $ is from Condition A3, and \begin{equation}
\mathcal{L}\left( \left( n^{-1/\alpha }l\left( n\right) \log \left(
\mathbf{Z}(nt),\mathbf{u}\right) ,0\leq t\leq 1\right)
\,|\,\mathbf{Z}(n)>0,\mathbf{Z}(0)=\mathbf{e}_{i}\right) \implies \mathcal{L}\left( L^{+}\right) .
\label{idva}
\end{equation}%
Combining (\ref{base}), (\ref{iras}), and (\ref{idva}) we obtain
(\ref{pria}) and (\ref{mash}). The theorem is proved.

\end{document}